\documentclass[12pt]{amsart}
\usepackage{amsfonts,newlfont,latexsym}
\usepackage{euscript}
\usepackage{amssymb,amscd}

      \theoremstyle{plain}
      \newtheorem{theorem}{Theorem}[section]
      \newtheorem{lemma}[theorem]{Lemma}
      \newtheorem{corollary}[theorem]{Corollary}
      \newtheorem{proposition}[theorem]{Proposition}
      \newtheorem{prop}[theorem]{Proposition}
      \newtheorem{remark}[theorem]{Remark}
      \newtheorem{conjecture}[theorem]{Conjecture}

\newtheorem{Theorem}{Theorem}[section]
\newtheorem{Proposition}[Theorem]{Proposition}
\theoremstyle{definition}
\newtheorem{definition}[theorem]{Definition}

\newcommand{\R}{\mathbb R}
\newcommand{\Rk}{\mathbb R^k}

\newcommand{\Z}{\mathbb Z}
\newcommand{\Zk}{\mathbb Z^k}

\newcommand{\Tk}{\mathbb T^k}

\def \a{\alpha}
\def \b{\beta}

\def \g{\gamma}
\def \D{\Delta}

\def \e{\varepsilon}

\def \A{(\EuScript{A})}

\def \dim{\mbox{dim}\,}
\def \d{\mbox{dist}\,}
\def \p{\mbox{Pr}\,}
\def \ker{\mbox{ker}}
\def \grad{\mbox{grad}\,}
\def \d{\mbox{dist}}

\def\PROOF{{\em Proof}\,: }
\def\Proof{{\em Proof}\,: }

\def\QED{~\hfill~ $\diamond$ \vspace{7mm}}

\def \A{\cal A}
\def \Rk {{\mathbb R}^k}
\def \rk {{\mathbb R}^k}

\def \rm {{\mathbb R} ^m}

\def \oa{{\cal O ^\ast}}
\def \o{{\cal O}}
\def \R{{\mathbb R}}
\def \Z{{\mathbb Z}}
\def \Zk{{\mathbb Z} ^k}

\def \la{{\lambda}}
\def \es{E ^s}

\def \w{{\cal W}}
\def \ws{{\cal W} ^s}

\def \eu{E ^u}

\def \wu{{\cal W} ^u }

\def \wi{{\cal W} _{\chi} }
\def \wh{{\cal W} _{H} }
\def \ci{C^{\infty}}

\begin{document}
\author[Boris Kalinin  and  Ralf Spatzier]
{Boris Kalinin $^\ast$ and Ralf Spatzier$^{\ast \ast}$}

\title[On the Classification of Cartan Actions]
{On the Classification of Cartan Actions}

\thanks{$^\ast$ Supported in part by NSF grants DMS-0140513}
\thanks{$^{\ast \ast}$ Supported in part by NSF grant DMS-0203735}

  \address{Department of Mathematics and Statistics,
  University of South Alabama, Mobile, AL36688}

  \email{kalinin@jaguar1.usouthal.edu}

 \address{Department of Mathematics, University of Michigan, Ann Arbor,
MI 48109.}

\email{spatzier@umich.edu}

\maketitle

\begin{abstract} We study higher rank Cartan actions on compact manifolds
preserving an ergodic measure with full support. In particular, we classify
actions by $\R ^k$ with $k \geq 3$ whose one-parameter groups act transitively  
as well as nondegenerate totally nonsymplectic $\Zk$-actions for $k \geq 3$. 
\end{abstract}

\section{Introduction}

The classification of Anosov systems is a deep and central problem in dynamics. 
For single diffeomorphisms, a long outstanding  conjecture asserts that they are 
all topologically conjugate to automorphisms of tori,  nilmanifolds and 
finite factors of such. 
Little progress has been made since Franks, Manning, and Newhouse proved this for Anosov 
diffeomorphisms on tori and nilmanifolds, and for codimension one Anosov 
diffeomorphisms \cite{Franks70,Manning73,Newhouse70}. There is no analogue 
to this  conjecture for flows. In fact, various 
examples of Anosov  flows  not topologically conjugate to an algebraic flow have 
been constructed by Franks and Williams \cite{FrW} and Handel and Thurston \cite{H-Th}.  
For single Anosov diffeomorphisms and flows one can easily change the derivative 
at periodic points. Thus topological conjugacies are rarely smooth. As Farrell and Jones 
have constructed Anosov diffeomorphisms on exotic tori, one cannot even hope for 
smooth classification of the underlying manifold structure \cite{FJ}.

The situation is quite different for higer rank Anosov actions, i.e. actions  of higher 
rank Abelian groups such that at least one element acts normally hyperbolically. 
The known actions enjoy very strong rigidity properties such as scarcity of invariant 
measures and cocycle rigidity (see surveys in \cite{Lind04,NT}). 
The easiest examples of such actions, arise
from products of Anosov diffeomorphisms or flows.  More interestingly, there are 
$\Zk$-actions on tori and nilmanifolds by automorphisms which are not products. 
There are also Anosov  $\Rk$-actions on homogeneous spaces $G/ \Lambda$  
by left translations, and more generally on biquotients. These form the class  of 
{\em algebraic actions}. Intriguingly, the only known examples of such actions 
are  either algebraic or are reducible, i.e. some finite cover admits an Anosov flow or 
diffeomorphism as a factor. By work of Palis and Yoccoz,  the centralizer of a generic 
Anosov diffeomorphism $f$ on a torus consists just of the powers $f^n$ of $f$ \cite{Palis-Yoccoz}.
A. Katok and the second author showed that  $C^1$-small 
perturbations of higher-rank algebraic Anosov actions with semisimple linear parts 
are {\em smoothly} conjugate to the original action \cite{KS97}. This followed 
earlier work by Katok and Lewis for the special case of a maximal commuting set 
of toral automorphisms \cite{KL1}.  Katok and Lewis also showed a global rigidity result for 
suitable higher rank actions on tori \cite{KL2}. Recently Damjanovic and Katok generalized
local rigidity to  partially hyperbolic actions on tori using KAM arguments \cite{DK}.  
Rodriguez Hertz classified Abelian actions  with an Anosov element for the 3-torus under
additional conditions on the action on homology.
These are all smoothly conjugate to a linear  action by automorphisms \cite{Hertz}.
All of these results motivate the following conjecture:

\begin{conjecture}  
All irreducible higher rank $\Zk$ and $\Rk$ Anosov actions on any compact manifold 
are {\em smoothly} conjugate to an algebraic action.

\end{conjecture}

While  this conjecture remains wide open in this generality, we will prove strong 
classification results for the subclass of Cartan actions in this paper. Cartan actions 
are Anosov actions such that the maximal non-trivial intersections of stable manifolds
of distinct elements are one-dimensional. This paper and part of its approach was 
motivated by similar results by E. Goetze and the second author  for Cartan actions 
by higher rank semisimple groups and their lattices \cite{GS}. 

Our main technical result is the following theorem. We call a one-parameter subgroup 
of a Lyapunov hyperplane in $\rk$ {\em generic} if it is not contained in any other 
Lyapunov hyperplane. Call an  $\Rk$ Cartan action {\em totally Cartan} if the set of Anosov
 elements is dense in $\Rk$.

\begin{theorem} \label{HolderMetric}
Let $\a$ be a 
totally Cartan action of $\Rk$, $k \geq 2$, on a compact 
smooth manifold $M$  preserving an ergodic probability measure $\mu$ with full 
support.  Suppose that every Lyapunov hyperplane contains a generic one-parameter
subgroup with a dense orbit.  Then there exists a H\"older continuous Riemannian 
metric $g$ on $M$ such that for any $a \in \rk$ and any Lyapunov exponent $\chi$
$$||a_\ast(v)||=e^{\chi (a)} ||v|| \qquad \text{ for any } \; v \in E_\chi.$$
\end{theorem}

We use this theorem to get the following classification   of $\rk$ Cartan actions actions.

\begin{theorem} \label{Main}
Let $\a$ be a $\ci$ totally Cartan action of $\Rk$, $k \geq 3$, on a compact
smooth connected manifold $M$ preserving an ergodic probability measure $\mu$ with full support.  
Suppose that every one-parameter subgroup of $\Rk$ has a dense orbit. 
Then $\a$ is $\ci$ conjugate to an almost algebraic action, i.e. the lift of the
 action to some finite cover of $M$ is $\ci$ conjugate to an $\Rk$-action by left
 translations on a homogeneous space $G/\Lambda$ for some Lie group  $G$ 
and cocompact lattice $\Lambda$. 
\end{theorem}

\medskip

The main new ingredient in the proof is the construction of the H\"{o}lder metric on 
the various Lyapunov foliations which is expanded and contracted precisely according to a 
linear functional. This is closely linked to cohomology triviality for cocycles. Indeed, 
Proposition \ref{HolderMetricSingle} says precisely that the restriction of the derivative 
cocycle in a Lyapunov
direction is  H\"{o}lder cohomologous to a linear functional. Cohomology triviality has 
been established for general cocycles for homogeneous actions \cite{KS95,KNT,NT}.
Nothing however seems to be known for general actions. Our approach here is specific
to the derivative cocycle, and is inspired by the proof of the Livsic' theorem.  Let us 
comment that in \cite{GS},  Goetze and the second author used topological super-rigidity 
techniques to trivialize the derivative cocycle for  actions of semi-simple groups. 
We also note that our results may provide another approach to the main results of 
\cite{GS}.

\medskip

Finally, we will apply  a technical variation Theorem \ref{technical} of our 
Theorem \ref {Main} to classify certain $\Zk$
Cartan actions. Call a Cartan action {\em totally nonsymplectic} or TNS if no two
nonzero Lyapunov exponents are negatively proportional. Further call an Anosov
action {\em non-degenerate} if the intersection of two Lyapunov hyperplanes is never 
contained in a third Lyapunov hyperplane. 

\begin{corollary} \label{MainZ}
Let $\a$ be a $\ci$ nondegenerate TNS Cartan action of $\Zk$, $k \geq 3$, 
on a compact smooth 
manifold $M$ such that each non-trivial element  is an  Anosov diffeomorphism. 
Suppose also that one of the diffeomorphisms is transitive. Then a finite cover of
$\a$ is $\ci$ conjugate to a $\Zk$ action by automorphisms of a nilmanifold.
\end{corollary}

In this corollary, the requirement that the action is TNS Cartan
is equivalent to the requirement that all Lyapunov exponents are simple and
there are no proportional Lyapunov exponents. 

The second author would like to thank E. Goetze for  numerous discussions related to this 
problem which yielded partial results and suggested  part of the approach in this paper.

\section{Basic Structures}

\subsection{Anosov actions}
Let us recall the definitions and basic properties of Anosov actions.

\begin{definition} {\sl Let $\alpha$ be a locally faithful action of $\Rk$  by 
smooth diffeomorphisms on a compact manifold $M$. 
Call an element $a \in \rk$  {\em Anosov} or {\em normally hyperbolic} for $\alpha$ 
if there exist real constants $\lambda>0$, $C>0$ and a continuous 
$\alpha$-invariant splitting of the tangent bundle
\[TM = E^u _a \oplus E^0  \oplus E^s _a\]
such that $E^0$ is the tangent distribution of the $\rk$-orbits, and for all $p \in M$, 
for all $v \in E^s _a (p)$  ($v \in E^u _a (p)$ respectively) and  $n >0$ ($n< 0$ 
respectively) the differential $a_* : TM \rightarrow TM$ satisfies
\[ \parallel a^n _* (v) \parallel \leq C e^{-\lambda  \mid n \mid } 
\parallel v 
\parallel. \]}
\end{definition}

Hirsch, Pugh and Shub introduced the notion of a diffeomorphism acting normally
hyperbolically with respect to an invariant foliation. Our Anosov elements are
precisely the elements in $\Rk$ which act normally hyperbolically with respect to
the orbit foliation of $\rk$ \cite{HPS}.  By \cite{HPS}, we can define stable  and 
unstable distributions $\es_a$ and $\eu _a$ for any  Anosov element $a \in \rk$.
These are H\"{o}lder distributions and integrate in the usual fashion to
stable and unstable foliations which we will denote by $\ws _a$ and $\wu _a$. 
These are H\"{o}lder foliations with $\ci$-leaves. 
(cf. \cite{HPS} for all this).

The set of Anosov elements ${\cal A}$ in $\rk$ is always an open subset 
of $\rk$ due to the structural stability theorem for normally hyperbolic maps 
by Hirsch, Pugh and Shub \cite{HPS}. 
 
\begin{definition}  Call $\alpha$ an {\em Anosov action} if some element 
$a \in \rk$ is Anosov. Furthermore call $\alpha$ {\em totally Anosov} if the 
set of Anosov elements $\cal A$ is dense in $\rk$. 
\end{definition}

It is not known if all Anosov actions are totally Anosov. We will assume 
henceforth that $\alpha$ is a totally Anosov $\rk$-action preserving
an ergodic probability measure $\mu$ with full support.

\subsection{Lyapunov theory}
 
Recall that for any diffeomorphism $\phi$ of a compact manifold $M$ preserving 
an ergodic  probability measure $\mu$, there are finitely many numbers $\chi ^i $ 
and a measurable splitting of the
tangent $TM = \bigoplus E^i$ such that the forward and backward Lyapunov
exponents of $v \in  E^i $ are exactly $\chi ^i $. This is the
Lyapunov decomposition of $TM$ for $\phi$. 

Now consider an $\rk$ action $\alpha$ on a compact
 manifold $M$ by diffeomorphisms preserving an ergodic probability 
measure $\mu$. Then we can refine the  Lyapunov decompositions
of the individual elements $ a \in \rk$ to a joint invariant splitting.

 \begin{prop} There are finitely many linear functionals $ \chi $ on $\rk$, a set of
 full measure ${\cal P}$ and 
  a measurable splitting of the tangent bundle $TM = \bigoplus E ^{\chi}$ over ${\cal
  P}$, invariant
  under $\alpha$,   such that 
  for all $a \in \rk$ and $ v \in E ^{\chi}$, the Lyapunov exponent of $v$ is $\chi 
  (a)$, i.e.
  \[ \lim _{n \rightarrow \stackrel{+}{_{-}} \infty } 
  \frac{1}{n} \log \parallel a ^n _\ast (v) \parallel = \chi (a)\]
  where $\parallel .. \parallel$ is some continuous norm on $TM$. 
\end{prop}
 
We call  $ \bigoplus E ^{\chi} $ the {\em Lyapunov splitting} and the nonzero 
linear functionals $\chi $ the {\em Lyapunov exponents} or {\em weights} of 
$\alpha$. We will call the hyperplanes $\ker \chi $ the {\em Lyapunov hyperplanes}
or {\em Weyl chamber walls}, and the connected components of $\rk
 - \cup _{\chi} \ker \chi $ the {\em Weyl chambers} of $\alpha$. 
 
Define the {\em coarse Lyapunov space} $E_{\chi} = \oplus E^{\lambda}$,
where the sum ranges over all positive multiples $\lambda = c \, \chi$ of $\chi$.
We will also denote $E_{\chi}$ by $E_H$ where  $H$ is the half space of $\rk$ 
on which $\chi$ is negative. Call such a half space $H$ a {\em Lyapunov half space}.
Note that $H$ is determined by the hyperplane $\ker \chi$ and an {orientation} of this 
hyperplane, i.e. a choice of one of the two half spaces  $\ker \chi$ bounds. In our case,
the orientation is given by which of the two half spaces $\chi$ is negative on. Then we 
obtain a measurable decomposition $TM = \oplus E_{H}$  where $H$ ranges over all 
Lyapunov half spaces. Note that for an $\R$-action we just retrieve the stable and 
unstable distributions.

\subsection{Coarse Lyapunov foliations} In this section we show that for a totally Anosov
$\rk$-action preserving an ergodic probability measure $\mu$ with full support the coarse 
Lyapunov splitting can be extended to a H\"{o}lder splitting of $TM$ consisting of
distributions tangent to foliations which we will call the coarse Lyapunov foliations. 
We also show that the Lyapunov hyperplanes, Weyl chambers, and
the coarse Lyapunov foliations agree for all invariant measures.

First we note that for each Anosov element $a \in \rk$ we have
$\es _a = \bigoplus_ {\chi (a) <0} E^{\chi}$ and $\eu _a = \bigoplus_ {\chi (a) >0} E^{\chi}$ 
at any point of the set ${\cal P}$ of full measure where the Lyapunov splitting is defined. 

\begin{proposition} \label{CoarseLyapunov}
Let $\alpha$ be a totally Anosov $\rk$-action preserving an ergodic  probability
measure $\mu$.
For each Lyapunov exponent $\chi $ and every $p \in {\cal P}$ 
\[E_{\chi}(p) = \bigcap _{\{a \in \A \mid \chi (a) <0\}} \es _a(p). \]
Moreover, if the measure $\mu$ has full support, the right hand side is H\"{o}lder 
continuous and thus $E^{\chi}$ can be extended to a H\"{o}lder distribution tangent 
to the foliation $\wi := \bigcap  _{\{a  \in \A  \mid \chi (a) <0\}} \ws _a$. In particular, 
$\wi $ has $\ci$-leaves.
\end{proposition}

\begin{remark}
{\em In this proposition the assumption that $\alpha$ is totally Anosov can be relaxed
to assuming the existence of an Anosov element in every Weyl chamber defined
by $\mu$.}

\end{remark}

\PROOF Let $\chi$ be a Lyapunov exponent, and define 
$L_ {\chi} = \bigcap _{\{a \in \A \mid \chi  (a) <0\}} \es _a $. This defines a distribution
on all of $M$. We will show that $L_ {\chi} = E{\chi}$ on ${\cal P}$ and that that $L_ {\chi}$
is continuous on the support of $\mu$, $supp \: \mu$.  Since $L_ {\chi}$ is the intersection 
of H\"{o}lder foliations, it follows easily that $L_ {\chi}$ is H\"{o}lder.

Let us first show that $L_ {\chi} = E{\chi}$ on the set of full measure ${\cal P}$ where the 
Lyapunov splitting is defined. Let  $p \in {\cal P}$. First  note that  
$E_ {\chi} (p) \subset \bigcap _{\{a \in \A \mid \chi  (a) <0\}} \es _a (p)$ since, by definition, 
$E_{\chi}$ is contained in every $E_a^s$ of the intersection.
To prove the reverse inclusion, suppose that for some $v \in T_p M$, 
$v \in \bigcap _{\{a \in \A \mid \chi (a) <0\}} \es _a (p)$ and $ v \notin E_{\chi} (p)$. 
Decompose $v = \sum _{\la }v_{\la}$ with $v_{\la} \in E_{\la} (p)$. Then $v_{\la} \neq 0$
for some $\la $ which is not a positive multiple of $ \chi $. Since $\A$ is dense in $\rk$, 
there is some $a \in \A$ such that $\chi (a) <0$ while $\la (a) >0$. This contradicts 
the fact that $ v \in \es _a$ and proves the reverse inclusion. 

Now we will show the continuity of $L_ {\chi}$ on $supp \: \mu$.
Since $\mu$ is ergodic, the dimensions of all the coarse Lyapunov distributions $E_ {\chi}$ 
and thus of all distributions $L_ {\chi}$ are constant on ${\cal P}$. It follows easily that these
distributions are continuous on ${\cal P}$ and form the direct sum 
$\bigoplus L_ {\chi}$ of dimension $(\dim M - k)$ over ${\cal P}$. Consider a point 
$q  \in  supp \:  \mu$ and a sequence of points $q_n \in {\cal P}$ converging to $q$. 
For each Lyapunov exponent $\chi$ let $R_{\chi}$ be a limit point of $L_{\chi} (q_n)$ in the
Grassman bundle of subspaces of dimension $\dim L_{\chi}$. To prove the continuity it
suffices to show that $R_{\chi} = L_{\chi} (q)$ for each $\chi$. By continuity of the $\es _a $, 
$R_{\chi} \subset L_{\chi} (q)$ for all $\chi$. To prove the reverse inclusion suppose that 
for some $\chi$ the dimension of $L_{\chi} (q)$ is greater than that of $R_{\chi}(q)$. 
Then the sum  $\bigoplus L_ {\chi}$ is no longer direct and for some $\chi$ the intersection 
$L_ {\chi} \cap \bigoplus _{\chi' \not= \chi} L_ {\chi '}$ is nontrivial. 
Let $v$ be a nonzero vector in this intersection and let $H$ be the negative half space of 
$\chi$.  Restrict the action to a generic 2-plane which intersects all Lyapunov half planes 
in distinct lines. Then it is easy to see that there are $a, b \in \A$ close
to $\partial H = \ker \la$,  such that $\chi (a) <0$ and $\chi (b) <0$ and such that $\chi$ 
and its positive multiples are the only Lyapunov exponents satisfying these two conditions. 
If we denote
\[ E_1= \es _a \cap \bigoplus _{\chi' \not= \chi} E_ {\chi '} \qquad \text{ and }
\qquad E_2=\eu _a \cap \bigoplus _{\chi' \not= \chi} E_ {\chi '} \]
these conditions imply that $\bigoplus _{\chi' \not= \chi} E_ {\chi '} = E_1 \oplus E_2$
and $E_1 \subset \eu _b$. We split $v$ as the direct sum $v=v_1+v_2$ and iterate it 
by $na$. Since $v \in L_\chi \subset \es _a$ and $v_1 \in E_1 \subset \es _a$ we see 
that $(na)v \to 0$ and $(na)v_1 \to 0$. Since $v_2 \in E_2 \subset \eu _a$ we conclude
that $v_2=0$ and $v=v_1$. Since $v \in L_\chi \subset \es _b$ and $v_1 \in E_1 \subset 
\eu _b$ we conclude that $v=v_1=0$. This shows that $R_{\chi} = L_{\chi} (q)$ for each 
$\chi$ and completes the proof of the proposition. 
\QED

\begin{lemma} \label{convexitylemma}
The set $\A$ of Anosov elements for $\alpha$ is the union of the
Weyl chambers in $\rk$.
\end{lemma}

\PROOF Suppose that $a, b \in \A$ belong to the same Weyl chamber $\cal C$.
Since $\es _a = \bigoplus_ {\chi (a) <0} E^{\chi}$, we get $\es _a = \es _b$. 
By commutativity, the stable distributions of Anosov elements are invariant under
$\alpha$. Let $c = t\, a + s \, b$ where $s, t >0$ are real numbers. 
If $ v \in \es _a (p)$ for $p \in M$, then the derivative $ c ^n _*  = (ns\, b) _* \circ  
(n t \, a)_*$ contracts $v$
exponentially fast as  $(n t \, a)_*$ does and 
$(n t \, a) _*(v) \in \es _b ((nt\, a)(p))$. We conclude that $\es _c$ is defined
and $\es _c = \es _a = \es _b$. Similarly, $\eu _c = \eu _a = \eu _b$ and hence
$c$ is Anosov. Thus the intersection of $\A$ with $\cal C$ is an open and dense 
convex cone in $\cal C$. Therefore $\cal C \subset \A$. Clearly, no element on a 
Weyl chamber wall can be in $\A$.  \QED

The elements of $\rk$ which belong to the union of the Weyl chambers are
called {\em regular}. All other elements of $\rk$ are called {\em singular}.
A singular element is called {\em generic} if belongs to only one Lyapunov
hyperplane. 

For an singular element $a \in \rk$ we can 
define its neutral, stable, and unstable distributions as
$$E_a^0 =T\o \oplus \bigoplus_ {\chi (a)=0} E^{\chi} \qquad
\es _a = \bigoplus_ {\chi (a) <0} E^{\chi} \qquad 
\eu _a = \bigoplus_ {\chi (a) >0} E^{\chi}.$$

\begin{lemma} \label{SingularElements}
Distributions $E_a^0$, $\es _a$, $\eu _a$ H\"older continuous. $\es _a$ and $\eu _a$
integrate to H\"older continuous foliations $\ws_a$ and $\wu_a$ with smooth leaves.
$\es _a$ is uniformly contracted and $\eu _a$ is uniformly expanded by $a$.
\end{lemma}
\begin{proof} 
The H\"older continuity of the distributions follows immediately from Proposition 
\ref{CoarseLyapunov}. To show the integrability we note that $\es _a = \bigcap \es_b$
and $\eu _a = \bigcap \eu_b$, where the intersection is taken over all Anosov
elements $b$ close to $a$. Indeed, if $b$ is close enough to $a$, the signs of 
$\chi (b)$ and $\chi (a)$ are the same for any Lyapunov exponent $\chi$ with 
$\chi (a) \not= 0$. This shows that $\es _a$ is contained in every $\es_b$ and
thus in the intersection.
For the reverse inclusion we note that for every nonzero Lyapunov exponent 
$\chi$ with $\chi(a)=0$ we can choose an Anosov element $b$ close to $a$ 
such that $\chi(b)>0$. The uniform contraction and  expansion can be obtained 
as in the proof of Lemma \ref{convexitylemma} since $a$ can be represented as
a positive combination of the nearby Anosov elements. \QED
 \end{proof}

\begin{remark} In contrast to the individual distributions $E^{\chi}$, $E_a^0$ is not 
necessarily integrable. Moreover, we do not assume any uniform estimates on the 
possible expansion or contraction of $E_a^0$ by $a$, so $a$ is not necessarily a 
partially hyperbolic element in the usual sense.
\end{remark}
\medskip

Now we will show that the structures of Lyapunov hyperplanes and Weyl chambers
agree for all invariant measures. Note that this does not entail that the Lyapunov 
functionals are the same. Indeed, they need not be as the case of products of Anosov
flows easily shows.

\begin{Proposition} \label{ConsistencyOfExponents}
Suppose that the Lyapunov splitting and the Lyapunov exponents exist
at a point $p$. Then the Lyapunov hyperplanes and Weyl chambers 
defined by these exponents coincide with the Lyapunov hyperplanes and 
Weyl chambers defined by the exponents of the ergodic invariant measure with
full support. Moreover, the coarse Lyapunov splitting of $T_pM$ defined by 
the exponents at $p$ coincides with the H\"older continuous coarse Lyapunov 
splitting defined in Proposition \ref{CoarseLyapunov}.
\end{Proposition}

\Proof  It suffices to show that for any Lyapunov half space $H$ defined by 
the ergodic invariant measure with full support and for any $v \in E_H(p)$ the 
Lyapunov exponent $\chi(\cdot, v): \rk \to \R$ has kernel $\partial H$ and
is negative on $H$. Suppose that this is not the case. Then there exists
$b \in H$ such that $\chi(b, v)>0$. Since the Anosov elements are dense
in $\rk$ we may choose $b$ to be Anosov, i.e. $b \in \A \cap H$. 
Then by the definition of $E_H$, $v \in E_H=\bigcap _{a \in \A \cap H} \es _a 
\subset \es _b$. But for $v \in \es _b$, $\chi(b, v)>0$ is impossible.
\QED

We immediately get that the Lyapunov half spaces on ${\cal P}$ are 
consistent with those  at all the periodic points.

\begin{corollary}
The Lyapunov hyperplanes, Weyl chambers, and coarse Lyapunov splitting 
for any compact orbit of the action coincide with the Lyapunov hyperplanes, 
Weyl chambers, and coarse Lyapunov splitting defined by the ergodic invariant 
measure with full support.
\end{corollary}

\subsection{Cartan  Actions}

Here we define  Cartan actions which are closely related to Hurder's trellised 
actions \cite{Hurder92,Hurder94}.

\begin{definition}
Call a (totally) Anosov action of $\rk$ a {\em (totally) Cartan} if
all   coarse Lyapunov foliations are one-dimensional. 
\end{definition}

Totally Cartan actions satsify the following properties  tantamount to being a 
trellised action.  First, let us call two foliations {\em pairwise transverse} if their 
tangent spaces intersect trivially.  This is different from the  standard notion in 
differential topology which also requires the sum of the tangent spaces to span 
the tangent space of the manifold. 

Consider a  totally Cartan action of $\Rk$ on $M$. Then the coarse Lyapunov 
foliations $\{\w_i \}$ form a collection of  one dimensional, pairwise transverse 
foliations such that

\begin{enumerate}
\item the tangent distributions have internal direct sum
$T\w_1 \oplus \cdots \oplus T\w_r \oplus T \o \cong TM$,
where $T\o$ is the distribution tangent to the $\rk$ orbits,
\item for each $x \in M$ the leaf $\w_i(x)$ of $\w_i$ through $x$
is a $C^\infty$ immersed submanifold of $M$,
\item the $C^\infty$ immersions $\w_i(x) \to M$ depend uniformly
H\"{o}lder continuously on the basepoint $x$ in the $C^\infty$ topology
on immersions, and
\item each $\w_i$ is invariant under every $a \in \rk$.
\end{enumerate}



\section{Proof of Theorem \ref{HolderMetric}} \label{ProofHolder}

It is clearly sufficient to show the existence of such a metric for each coarse
Lyapunov distribution. After that, the desired metric on $M$ can be obtained
using these metrics and the natural metric on the orbit distribution, by
requiring that the coarse Lyapunov distributions and the orbit distribution
are pairwise orthogonal. The metric on the orbit distribution can be defined
as follows. For $v \in T_x \o$ define $||v||=||b||_{\Rk}$, where $b \in \Rk$ 
is such that $v=\frac{d}{dt}((tb)x)$.
Thus the theorem reduces to the following proposition.

\begin{proposition} \label{HolderMetricSingle} 
Let $\chi$ be a Lyapunov exponent, $H$ be its negative Lyapunov half space,
and $E=E_H$ be the corresponding one-dimensional coarse Lyapunov distribution. 
Under the assumptions of  Theorem \ref{HolderMetric} there exists a H\"older 
continuous Riemannian metric on $E$ for which 
$$||a_\ast(v)||=e^{\chi (a)} ||v||$$ 
for any $a \in \rk$ and $v \in E$.
\end{proposition}

\begin{proof}
Let $\w=\wh$ be the coarse Lyapunov foliation of the coarse 
Lyapunov distribution $E$. Denote by $E'$  and $\w'$  the (possibly trivial) coarse 
Lyapunov distribution and the coarse Lyapunov foliation corresponding to the
Lyapunov half space $-H$.

{\bf Notations.} In Sections \ref{ProofHolder} and \ref{ProofClosing} for any element 
$b \in \Rk$ we denote by $D_x^E b$ the restriction of its derivative at $x\in M$ to 
$E(x)$.  We fix some background Riemannian metric $g_0$  and denote the norm 
of $D_x^E b $ with respect to $g_0$ by $d^E_x b = || D_x b (v) ||_{bx}  \cdot ||v||_x^{-1}$, 
where $v \in E_x$ and $||.||_x$ is the norm given by $g_0$ at $x$.

By the assumption there exist an element $a_0 \in  \partial H$ not contained in any other
Lyapunov hyperplane and a point $x^\ast$ such that the orbit $\oa :=\{(ta_0)x^\ast \}$ is dense.
We define a new metric $g^\ast$ on $E$ over $\oa$ by taking the background metric 
$g_0$ on $E_{x^\ast}$ and propagating it along this dense orbit by the derivative 
$D^E_{x^\ast}(t a_0)$. By the construction, the derivative $D_x^E(t a_0)$ is isometric 
with respect to $g^\ast$ for any $t$ and any $x\in \oa$. 

The main part of the proof is to show that the metric $g^\ast$ is H\"older continuous 
on $\oa$ and thus extends to a H\"older continuous Riemannian metric $g$ on the whole
distribution $E$. Clearly, such $g$ is also preserved by $D_x^E(t a_0)$ for any $t$ and 
any $x\in M$. For any other element $b \in  \rk$ consider the metric $b_\ast g$. 
By commutativity, this metric is again preserved by $ta_0$ for any $t$. Since $E$ is 
one-dimensional and $\oa$ is dense, it is easy to see that, up to constant scaling, there 
is only one metric invariant under $ta_0$. Hence $b_\ast g = c \cdot g$ on $M$, where 
$c$ is a costant. Clearly, the logarithm this constant gives the Lyapunov exponent for $b$
of any $v \in E$ and hence $c=e^{\chi (b)}$. To complete the proof of the proposition we 
will now show that $g^\ast$ is H\"older continuous on $\oa$.
\medskip

To prove that $g^\ast$ is H\"older continuous on $\oa$ we need to show that for any point 
$x \in \oa$ which returns close to itself under an element $a=ta_0$ the norm $d^E_{x} a$
defined above is H\"older close to $1$. The specific H\"older exponent depends on the 
H\"older exponents of certain invariant foliations. Let $\a_0>0$ be such that all coarse 
Lyapunov distributions are H\"older continuous  with exponent $\a_0$, and let 
$\a=\min \{\a_0, \frac12\}$. We will show that for any positive $\b <\a/(1+\a)$ there exists 
a positive constant $\e_\ast$ such that  
\begin{equation}\label{Holder1}
|d^E_x a-1| < \d (x, ax)^\b \qquad 
\text{ for any } x \in \oa \text{ and } a \text{ with } \d (x, ax)<\e_\ast
 \end{equation}

To show this we will use special  closing arguments given in 
Propositions \ref{Closing0} and \ref{Closing1}. They establish the existence 
of a nearby point which returns under $a$ to the same leaf of $\o \oplus \w'$.

Fix any positive $\b <\a/(1+\a)$. Then $1-\b>1/(1+\a)$. 
Fix a positive $\g$  smaller than $1-\b$ but  greater than $1/(1+\a)$. From this we obtain 
\begin{equation}\label{Holder1''}
\a \g > \a/(1+\a)>\b
 \end{equation}
 
For such $\b$ and $\g$, Proposition \ref{Closing0} gives a
positive constant $\e_0$.  We choose $\e_\ast>0$ such that $\e_\ast < \e_0$
and $\e_\ast^{\gamma} < \e_1$, where $\e_1$ is a positive constant given
by Proposition \ref{Closing1}.

Consider $x_0 \in \oa$ and $a=ta_0$ with $\e=\d (x_0, ax_0)<\e_\ast$.
It suffices to assume that the return time $t$ is positive and large enough.
Suppose that the inequality \eqref{Holder1} does not hold for $x_0$.
We will use Propositions \ref{Closing0} and \ref{Closing1} to obtain an 
estimate for $d^E_{x_0} a$ which contradicts this assumption if $\e_\ast$
satisfies the inequality \eqref{Holder8}. This will imply that with such $\e_\ast$
the inequality \eqref{Holder1} holds for any 
$x \in \oa$ and $a=t a_0$ with $\d (x, ax)<\e_\ast$.
\medskip

If the inequality \eqref{Holder1} does not hold for $x_0$, by taking inverses 
we may assume without loss of generality that $d^E_{x_0} a <1-\d (x_0, ax_0)^\b$.
Thus we can use Proposition \ref{Closing0} to obtain the corresponding point $x_1$. 
Since $\d (x_1 , a x_1)<\e_\ast^{\gamma} < \e_1$ by the choice of $\e_\ast$, 
we can use Proposition \ref{Closing1} to obtain the corresponding point $x_2$
and $\delta \in \Rk$ for which $(a+\delta)x_2 \in \w'(x_2)$. 

Denote $b=a+\delta$. Take an Anosov element $c$ such that $\ws_c=\ws_a \oplus \w'$. 
Let $y=\lim (t_n c) x_2$ be an accumulation point of $c$-orbit of $x_2$.
Since $b x_2  \in \w'(x_2)$ and $c$ contracts $\w'$ we obtain $y=\lim (t_nc)(bx_2)$.
Then $by=\lim b((t_nc) x_2)=\lim (t_nc )(bx_2)=y$ and thus $y$ is a fixed point for $b$. 

We will now show that $d^E_y b$ is close to $1$. 
Since $b-a=\delta$ is small, we may assume that $b$ does not belong 
to any Lyapunov hyperplane different from $\partial H$. Hence either $b$ is Anosov
or $b$ belongs to $\partial H$. In the latter case, since $by=y$ we immediately obtain
$d^E_y b=1$. Indeed, suppose for example that $d^E_y b=\lambda>1$. Arbitrarily close to
$b$ there are Anosov elements for which $E$ is contained in the stable distribution.
For any such element $c$ and for $n$ sufficiently large so that $-nc$ expands $E$ 
we can estimate $d^E_y (n(b-c))= d^E_y (-nc) \cdot d^E_y (nb) \ge \lambda ^n$, 
but this is impossible if the element $b-c$ is sufficiently close to $0\in \rk$.
If, on the other hand, $b$ is Anosov we can conclude that the orbit $\Rk b$ is compact 
(cf. \cite{Qian94}). Hence the Lyapunov splitting and Lyapunov exponents are defined
at all points of this compact orbit. We denote by $\tilde \chi$ the Lyapunov exponent of
vectors in $E$ on this compact orbit. Note that while $\tilde \chi$ may not coincide with
$\chi$, by Proposition \ref{ConsistencyOfExponents} their kernels are the same:
 $\ker \tilde \chi =\partial H$. Since $y$ is fixed by $b$, 
$\tilde \chi (b) =\log (d^E_y b)$. Since $a \in \ker \tilde \chi =\partial H$, we obtain 
$$|\tilde \chi (b)| =|\tilde \chi (a) + \tilde \chi (\delta)|= |\tilde \chi (\delta)| < 
C_1 ||\delta|| < C_2  \e^{\gamma} \text{ , where } \e=\d (x_0, ax_0)$$
Thus we conclude that
\begin{equation}\label{Holder2}
|d^E_{y} b-1| < C_3 \e^{\gamma}
 \end{equation}
 
We will now show that $d^E_{x_0} a$ is H\"older close to $d^E_y b$.
using the following estimates.
First we get from parts (1) and (2) of Proposition \ref{Closing0} together
with Lemma \ref{d^Ea HolderContinuity} that
\begin{equation}\label{Holder3}
|d^E_{x_0} a-d^E_{x_1} a| < C_4 \e^{\a\gamma} 
 \end{equation}
Next combine  parts (1) and (2) of Proposition \ref{Closing1} together
with the analog of Lemma \ref{d^Ea HolderContinuity} for $\wu _{a_0}$.
\begin{equation}\label{Holder4}
|d^E_{x_1} a-d^E_{x_2} a| < C_5 \e^{\a\gamma} 
 \end{equation}
Now apply part (5) of Proposition \ref{Closing1} with $b-a=\delta$.
\begin{equation}\label{Holder5}
|d^E_{x_2} a-d^E_{x_2} b| < C_6 \e^{\gamma} 
 \end{equation}
 
 Finally, we show that
 \begin{equation}\label{Holder6}
|d^E_{x_2} b-d^E_{y} b| < C_7 \e^{\gamma} 
 \end{equation}
This can be seen as follows. By commutativity $b=(-t_nc) \circ b \circ (t_nc)$,
so for any iterate $(t_nc) x_2$ we have
$$d^E_{x_2} b= d^E_{b(t_nc)x_2} (-t_nc)   \cdot d^E_{(t_nc)x_2} b  \cdot d^E_{x_2} (t_nc)$$ 
The middle term on the right side tends to $d^E_{y} b$, while the ratio of the
other two terms is H\"older close to $1$. The latter follows from a standard argument in
the proof of Livsic' theorem since $d^E_x c$ is H\"older continuous and the orbits 
of $x_2$ and $b x_2$ under $c$ are exponentially close.

We conclude that equations \eqref{Holder2},  \eqref{Holder3},  \eqref{Holder4},  
\eqref{Holder5},  \eqref{Holder6} imply 
\begin{equation}\label{Holder7}
|d^E_{x_0} a-1| < C_8 \e^{\a\gamma} 
\end{equation}

However, by equation \eqref{Holder1''},  $\a \gamma  >  \b$. Thus, possibly 
after decreasing $\e_\ast$ further to satisfy
\begin{equation}\label{Holder8}
C_8 \e_\ast^{\a\gamma} < \e_\ast^\b
\end{equation}
we obtain that the inequality \eqref{Holder7} contradicts the assumption that
the inequality \eqref{Holder1} does not hold for $x_0$. Hence we conclude that 
with this $\e_\ast$ equation \eqref{Holder1} holds for any $x \in \oa$.

This establishes the desired H\"older estimate for the constructed metric and 
thus completes the proof of Proposition  \ref{HolderMetricSingle} and 
Theorem \ref{HolderMetric}.
\QED

\end{proof}

\section{Closing lemmas}  \label{ProofClosing}
In this section we fix a coarse Lyapunov subspace $H$ and a singular element 
$a_0 \in  \partial H$ which is generic, i.e. is not contained in any other Lyapunov
hyperplane. Hence the neutral distribution of $a_0$ is $E_H \oplus E_{(-H)}$,
where $E_{(-H)}$ is trivial if there is no Lyapunov exponent positive on $H$.
We will use notations $E=E_H$, $E'=E_{(-H)}$, $E^s=E^s_{a_0}$,  and 
$E^u=E^u_{a_0}$. We assume that the corresponding distributions $\w$, $\w'$, 
$\ws$, and $\wu$ are H\"older continuous with exponent $\a_0>0$, and denote 
$\a=\min\{\a_0, \frac12\}$.

The propositions below can be compared to the Anosov closing lemma for Anosov flows.
The main differences are the following. The neutral distribution $E \oplus E'$ of
$a_0$ is not integrable in general, this forces us to consider $\w$ and $\w'$
separately. The holonomies under consideration are only H\"older continuous,
this forces us to use a topological fixed point argument rather than the contracting
mapping theorem. Proposition \ref{Closing1} can be formulated and proved in the 
context of a single partially hyperbolic diffeomorphism $a_0$ to become a relatively 
standard version of the closing lemma for diffeomorphisms normally hyperbolic to 
an invariant foliation (cf. \cite{HPS}). We formulate it in the specific form that 
we need with the Lipschitz estimates absent in \cite{HPS}. The main novelty is
in Proposition \ref{Closing0} where we utilize the weak contraction in the neutral
foliation $\w$. This is substantially higher rank, since the proof relies on the nonstationary
linearization along the leaves of $\w$ used in Lemma \ref{d^Ea HolderContinuity}. 
The existence of this nonstationary linearization is provided by an element in $\rk$
which uniformly contracts $\w$.

\begin{proposition} \label{Closing0}
For any positive $\b <\a/(1+\a)$ and positive $\g<1-\b$ there exists $\e_0>0$ 
such that for any $x_0 \in M$ and $t >1$ with 
$\e = \d (x_0 , (ta_0) x_0)<\e_0$ and $||D^E_{x_0} (ta_0)||<1-\e^\b$
there exists a point $x_1 \in M$  such that 

\begin{enumerate}
\item $\d (x_0,x_1) <  \e^{\gamma}$
\item $x_1 \in (\ws \oplus \w)(x_0)$
\item $\d ((ta_0) x_0 , (ta_0) x_1)< \e^{\gamma}$
\item $(ta_0)x_1 \in (\o \oplus \wu \oplus \w')(x_1)$
\item $\d (x_1 , (ta_0)x_1)< \e^{\gamma}$
\end{enumerate}
\end{proposition}

\begin{proposition} \label{Closing1}

There exist positive constants $\e_1$ and $C$ such that  
for any $x_1 \in M$ and $t>1$ with 
$(ta_0) x_1 \in (\o \oplus \wu \oplus \w')(x_1)$ and $\e=\d (x_1 , (ta_0) x_1)<\e_1$, 
there exist a point $x_2 \in M$, $\delta \in \Rk$, such that 

\begin{enumerate}
\item $\d (x_1,x_2) <  C\e$
\item $x_2 \in \wu (x_1)$
\item $\d ((ta_0)x_1,(ta_0)x_2) <  C\e$
\item $(ta_0+\delta)x_2 \in \w'(x_2)$
\item $||\delta|| < C\e$
\end{enumerate}
\end{proposition}

Below we give a proof of Proposition \ref{Closing0}.  Proposition \ref{Closing1} can 
be proved similarly, its proof would avoid the main difficulty caused by the neutral
direction $\w$ and yield the Lipschitz estimates.  

\medskip

\begin{proof}
Recall that 
$\b <\a/(1+\a)$ implies $1-\b>1/(1+\a)$. 
Since it is clearly sufficient to consider only $\g$ which are close to $1-\b$,
we may assume that $\g>1/(1+\a)$. From this we obtain that 
$\g>1-\a$ and $\a \g > \a/(1+\a)>\b$.
 We summarize the  inequalities we have:
 \begin{equation}\label{closing1}
 0<\b<\frac{\a}{1+\a}<\a\g < \a \le \frac12 \le 1-\a < \g < 1-\b <1
 \end{equation}

We introduce the following notations $a=ta_0$, $y_0=a x_0$, and $F=\ws \oplus \w$.
By assumption $\e = \d (x_0 , y_0)<\e_0$. Consider balls $B_1 \subset \ws (x_0)$ and 
$B_2 \subset \w (x_0)$ of radius $k_1 \e^\g$ centered at $x_0$. For $w_1 \in B_1$ and 
$w_2 \in B_2$ we denote by $[w_1,w_2]$ the unique local intersection of $\w(w_1)$ and
$\ws (w_2)$ in $F=\ws \oplus \w$.  We define a "rectangle"
$P=\{[w_1,w_2] :  w_1\in B_1,w_2\in B_2\} \subset F(x_0)$. 
Since the minimal angle between the leaves of foliations $\w$ and $\ws$ is bounded 
away from $0$, the constant $k_1$ can be chosen so small that  $P$ is contained in 
a ball of radius $\e^\g$ centered at $x_0$. Denote $f=a |_P : P \to F(y_0)$ 
and let $h  : f(P) \to F(x_0)$ be the holonomy map of the foliation 
$\o \oplus \wu \oplus \w'$. We will show that for small enough $\e_0$ we can ensure
that  $h(f(P)) \subset P$. Since $P$ is homeomorphic to a ball, this implies the existence 
of a fixed point  $x_1$ for $h \circ f$ which satisfies the conclusions of the proposition.  
\medskip

We identify a neighborhood of $x_0$ in $M$ with $T_{x_0} M$ using
local coordinates for which the differential at the base point is identity,
and the leaf $F(x_0)$ identifies with its tangent space.  
By abuse of notations we will write $F$ for the flat leaf $F(x_0)$ and $F'$ for
the leaf $F(y_0)$. We denote by $\p$ the orthogonal projection from $F'$ to $F$.
\medskip

First we give an estimate of the distance between the leaves $F$ and $F'$. For a
point $x \in F$ let $x' \in F'$ be such that $x=\p (x')$. Denote by $d(x)$ the distance
between $x \in F$ and the corresponding $x' \in F'$. Let $d_0= d(x_0)$. Since the 
tangent distribution of the foliation $F$ is H\"older with exponent $\a$, the maximal 
angle between $F$ and $T_{x'}F'$ is at most $K_2 d(x)^\a$. Hence the function $d$
satisfies inequality $|\grad d(x)|\le K _2d(x)^\a$. We conclude that $d$ is bounded 
by the solution of $y'=K_2 y^\a$ with the initial condition $y(0)=d_0\;$:
$$ d(x)\le \left( d_0^{1-\a}+(1-\a)K_2 \d (x, x_0) \right) ^{1/(1-\a)}.$$

We observe that $d_0$ is of order $\e = \d (x_0 , y_0)$. Recall that $P$ is contained 
in a ball of radius $ \e^\g$ centered at $x_0$, and that $\e^\g<\e^{1-\a}$ 
by \eqref{closing1}. Thus the the first term in the sum dominates for small $\e$ and
hence the maximum of $d$ on $P$ is bounded by $K_3 \e$ for some constant $K_3$.
We conclude that 
\begin{equation}\label{closing1.5} 
\d (P,f(P)) = \max  \d (x, x')  \le K_3 \e.
 \end{equation}
\medskip

Now we estimate how the holonomy $h$ from $F'$ to $F$ deviates from the 
orthogonal projection $\p$. Since the minimal angle between the leaves of
foliations $F$ and $\o \oplus \wu \oplus \w'$ is bounded away from $0$, 
there exists a constant $K_4$ such that for any 
$x' \in f(P)$ we can estimate 
\begin{equation}\label{closing2}
\D(x') :=\d (h(x'), \p (x'))= \d (h(x'), x))\le K_4 \, \d(x,x') \le K_4  K_3 \e 
 \end{equation}
\medskip

Now we will study $f(P)$ and its projection $\p (f(P))$ to $F$. Our goal is to show
that $\p (f(P)) \subset P$ and that the distance form $\p (f(P))$ to the relative 
boundary $\partial P$  of $P$ in $F$ is greater than the upper bound we just 
obtained for the distance $\D$ between the holonomy $h$ and the projection 
$\p$. This will imply that $h(f(P))\subset P$ and complete the proof.

We first estimate the derivatives of $a=ta_0$ on $P$. Recall that $a_0$ contracts 
$\ws$. We may assume that the return time $t$ has to be large for the return under 
$a$ to be $\e$-close. Hence we may assume that the norm of the derivative of $a$ 
restricted to $E^s$ is bounded above by $\frac14$ on $P$ : 
\begin{equation}\label{closing2.5}
||D^{E^s}_x a|| < \frac 14 \qquad \text{for any } x\in P 
 \end{equation}

Now we will estimate the derivative $D_x^E a$ in $E$ direction using the following lemma.

\begin{lemma} \label{d^Ea HolderContinuity}
The dependence of the derivative $D^E_y (ta_0)$ on $y$
is Lipschitz continuous along $\w$ and H\"older continuous 
along $\ws _{a_0} $ with exponent $\a$ and constants 
independent of $y$ and $t$.
\end{lemma}
\begin{proof}
To show the Lipschitz continuity along one-dimensional leaves of $\w$ we use 
nonstationary linearization of the action. For an element $b \in \Rk$ which 
contracts $\w$ the following lemma from \cite{KL2} gives the nonstationary
linearization of $b$ along $\w$.

\begin{lemma} \label{Linearization}
If a diffeomorphism $b$ of a manifold $M$ contracts an invariant
one-dimensional foliation $\w$, then there exists a unique family of 
$C^\infty$ diffeomorphisms $h_x: \w(x) \to T_x\w$, $x\in M$, 
such that

 (i) $\;\;\; h_{b x}\circ b=D_xb  \circ h_x$,
 
 (ii) $\;\; h_x(x)=0$  and $D_xh_x$ is the identity map, 
  
 (iii) $\; h_x$ depends continuously on $x$ in $C^\infty$ topology.
\end{lemma}

Since $a=ta_0$ commutes with $b$, the same family linearizes $a$ and we obtain
$$a |_{\w(x)} (y) = (h_{a x} ^{-1} \circ D^E_x a  \circ h_x )(y) : \w(x) \to \w(ax)$$
Since the second derivatives of $h_x$ are uniformly bounded, the first
derivatives vary Lipschitz continuously. This implies that $D^E_y a$ varies
Lipschitz continuously along $\w(x)$ in a small neighborhood of an
arbitrary point $x$.

Now we show the H\"older continuity along $\ws$.
Since $a_0$ exponentially contracts $\ws_{a_0}$, the orbits under $t_0a$ of any 
two nearby points $y_1$ and $y_2 \in \ws_{a_0} (y_1)$ are exponentially close. 
The derivative cocycle in the direction of $E$ is a H\"{o}lder cocycle Hence the 
standard argument from the proof of Livsic' theorem shows that $D^E_{y_1}(ta_0)$
is H\"older close to $D^E_{y_2}(ta_0)$ with exponent $\a$ and a uniform constant.
 \QED
\end{proof}

By the assumption, $||D^E_{x_0} a|| <1-\e^\b$. Since $P$ is contained in
a ball of radius $\e^\g$, using Lemma ~\ref{d^Ea HolderContinuity} and the 
fact that $\b<\a \g$ we obtain
\begin{equation}\label{closing3}
||D^E_{x} a|| < 1- \frac{\e^\b}2 \qquad \text{for any } x\in P 
 \end{equation}
provided that $\e_0$ is small enough. 

\medskip
 
Now we study how $f(P)$ projects to $F$ and estimate the distance form 
$\p (f(P))$ to the boundary $\partial P$. 
Recall that $B_2$ is just an interval in $\w(x_0)$ centered at $x_0$ of length 
$2k_1\e^\g$. Using estimates similar to \eqref{closing1.5} we obtain 
\begin{equation}\label{closing4} 
\d (B_2,\p (f(B_2)) \le \d(B_2,f(B_2)) \le K_3 \e.
 \end{equation}
The "rectangle" $P$ is the union of "horizontal layers" 
$B_x=B_1 \times \{x\} =\{[w_1,x]_b :  w_1\in B_1\} \subset \ws (x)$ for $x \in B_2$.
Let us fix $x \in B_2$ and the corresponding $B_x$.  Using \eqref{closing2.5} we
obtain that $f(B_x)$ is contained in a ball of radius $\frac13 k_1\e^\g$ centered
at $y=f(x)$ in $\ws (y)$. By \eqref{closing4}, $\d(\p (y),B_2) \le K_3 \e$. Hence
it is easy to see that the projection of this ball to $P$ is at least $k_5 \, \e^\g$
away from the "vertical part" $\partial B_1 \times B_2$ of the boundary 
$\partial P$. Thus we obtain
\begin{equation}\label{closing5} 
\d (\p (f(B_x), \partial B_1 \times B_2) \ge k_5 \, \e^\g.
 \end{equation}

Now we will estimate the distance from $\p (f(B_x))$ to $B_{x_1} =
B_1 \times \{x_1\} \subset \partial P$, where $x_1$ is one of the endpoints of 
the interval $B_2$. We denote $z=\p (y)$, and $z_0=\p (y_0)$. The derivative 
estimate \eqref{closing3} above implies that $\d (y_0,y) \le k_1\e^\g(1- \frac12 \e^\b)$, 
hence after projecting we have $\d (z_0,z) \le k_1\e^\g(1- \frac12 \e^\b)$. Since 
$\d(x_0,y_0)\le \e$ by assumption, we obtain $\d(x_0,z_0)\le \e$. 
Hence $\d (x_0,z) \le \e+k_1\e^\g(1- \frac12 \e^\b)$ and 
\[ \d(z,x_1) \ge \d (x_0,x_1)-\d (x_0,z)\ge k_1\e^\g -\left(\e + k_1\e^\g - 
\frac12 k_1\e^\g \e^\b  \right) \ge \frac12  k_1\e^\g \e^\b -\e \ge \frac13 k_1\e^{\g+\b}  \]
since $\g+\b<1$, provided that $\e<\e_0$ is small enough.  By \eqref{closing4},
$z$ is $K_3 \e$ close to the interval $B_2 \subset \w (x_0)$. Since the angles 
between $\w$ and $\ws$ are bounded away from zero, it is easy to see that
\[ \d(z,B_{x_1}) \ge k_6 \e^{\g+\b} \]

We will now complete the estimate of the distance between $\p (f(B_x))$ and 
$B_{x_1}$. Note that $f(B_x)$ and $B_{x_1}$ lie on two nearby leaves of the 
foliation $\ws$ and the distance from any point in $\p (f(B_x))$ to any point
$B_{x_1}$ is at most of order $\e^\g$. By H\"older continuity of the corresponding 
distribution we see that the angles between tangent spaces to $f(B_x)$ and 
$B_{x_1}$ differ no more than $K_7 (\e^\g)^\a$. Hence, after projecting, the 
angles between tangent spaces to $\p (f(B_x))$ and $B_{x_1}$ also differ 
no more than $K_7 (\e^\g)^\a$. Thus, the distance from a point on $\p (f(B_x))$
to $B_{x_1}$, as a function of this point, cannot change by more than 
$$K_7 (\e^\g)^\a \cdot (2k_1\e^\g ) = 2k_1 K_7  \e^{\g+\a\g}.$$
Since $\a\g>\b$ we obtain 
\[ \d(\p ((f(B_x)), B_{x_1}) \ge  k_6 \e^{\g+\b} - 2k_1K_7  \e^{\g+\a\g} \ge k_8 \e^{\g+\b} \]
provided that $\e<\e_0$ is small enough. Combining this with \eqref{closing5} 
we conclude that 
$$
\d (\p (f(B_x), \partial P) \ge \max \{ k_5 \, \e^\g, k_8 \e^{\g+\b} \}.
$$
Since $x \in B_2$ was arbitrary we conclude that 
 $$\d (\p (f(P), \partial P) \ge  k_9 \e^{\g+\b} .$$

Since $\g+\b<1$ we see this distance is larger than that the estimate \eqref{closing2}
for the deviation of the holonomy from the projection, provided that $\e<\e_0$ is small 
enough. This shows that $h(f(P))\subset P$ and proves the existence of a fixed point 
which satisfies the conclusions of Proposition  \ref{Closing0}. \QED
\end{proof}



\section{Classification of $\rk$-actions}

In this section we will prove the following generalization of Theorem \ref{Main}.
We will use this stronger technical version to prove Corollary \ref{MainZ} in the next section.
Given two Lyapunov hyperplanes $H_i$ and $H_j$ with $H_j \not= \pm H_i$, we denote by 
$E_{ij}$ the smallest direct sum of coarse Lyapunov distributions which is integrable and 
contains $E_{H_i}$ and $E_{H_j}$. We denote the foliation tangent to $E_{ij}$ by $\w _{ij}$

\begin{Theorem}\label{technical}
Let $\a$ be a $\ci$ totally Cartan action of $\Rk$, $k \geq 3$, on a compact smooth 
connected manifold $M$ preserving an ergodic probability measure $\mu$ with full support.  
Suppose that every Lyapunov hyperplane contains a generic one-parameter subgroup with 
a dense orbit.  Assume further that for any two Lyapunov half spaces $H_i$ and 
$H_j \not= \pm H_i$ there is an element $a \in \partial H_i \cap \partial H_j$ and a point 
$x \in M$ for which the closure of the $a$-orbit contains the whole leaf $\w_{ij} (x)$. 
Then $\a$ is $\ci$ conjugate to an almost algebraic action.
\end{Theorem} 

\begin{proof} First we note that Theorem \ref{HolderMetric} provides us with a H\"older 
continuous Riemannian metric $g$ on $M$ such that for any $a \in \rk$ and any Lyapunov 
exponent $\chi$
$$||a_\ast(v)||=e^{\chi (a)} ||v|| \qquad \text{ for any } \; v \in E_\chi.$$

The next step is to show that this metric $g$ and the coarse Lyapunov splitting are $\ci$. 
For this we will use Theorem 2.4 from \cite{GS}. Our assumptions on the one-dimensionality 
of the coarse Lyapunov foliations and the existence of the metric $g$ given by Theorem
\ref{HolderMetric} guarantee that the assumptions of Theorem 2.4 are satisfied with the 
exception of the assumptions of invariant volume and ergodicity of one-parameter subgroups. 
We note however that the proof goes through verbatim under the weaker assumption of preservation
of an ergodic invariant measure with full support.  Furthermore, ergodicity of one-parameter 
subgroups is used only to guarantee that for certain one-parameter subgroups in the Lyapunov
hyperplanes there are points $x$ whose orbits accumulate on the whole leaves $ \w _{ij} (x)$.
The last assumption of Theorem \ref{technical} is precisely what is required in the proof.
Thus we obtain that the metric $g$ and the coarse Lyapunov splitting are $\ci$. Now we can
complete the proof of Theorem \ref{technical} as follows.

Passing to a finite cover
of $M$ if necessary, we may assume for any coarse Lyapunov direction $E_{\chi}$
that there are nowhere vanishing  vectorfields tangent to $E_{\chi}$.
 Consider all the vectorfields
$V_{\chi}$  pointing in the various one-dimensional Lyapunov foliations ${\cal W}_{\chi}$
of length 1 with respect to the metric $g$ as well as the generating fields of the 
$\Rk$-action. Then all vectorfields $V_{\chi}$ are smooth and are expanded or 
contracted uniformly by $e^{\chi (a)}$  for $a \in \Rk$. Hence the Lie bracket of 
two such fields $[V_{\chi}, V_{\xi}]$ is expanded by $e^{(\chi + \xi)(a)}$, and hence 
is a constant multiple of $V_{\chi + \xi}$ if the latter exists. Otherwise 
$ [V_{\chi}, V_{\xi}]=0$. Since the $\Rk$-action normalizes the $V_{\chi}$, 
the $V_{\chi}$ together with the the generating fields of the $\Rk$-action span a 
finite dimensional Lie algebra $\mathfrak{g}$. Let $G$ be the cprresponding simply connected 
Lie group 
with Lie algebra $\mathfrak{g}$. Since the $V_{\chi}$ are globally defined and bounded 
with respect to an ambient Riemannian metric, $G$ acts on $M$. By the construction 
this action is locally simply transitive and thus transitive as $M$ is connected. Hence $M$ is a homogeneous 
space $G/ \Gamma$ for a lattice $\Gamma$ in $G$. Moreover, the $\Rk$-action embeds 
into the left action of $G$ by construction. \QED

\end{proof}



\section{Proof of Corollary \ref{MainZ}}

To prove Corollary \ref{MainZ}  we will apply Theorem \ref{technical} to the suspension
$\a '$ of the $\Zk$ action $\a$. By the assumption $\Zk$ contains an element $a$ which
is a transitive Anosov diffeomorphism. It is well known that such an element 
has a
unique measure $\mu$ of maximal entropy. In fact, transitive Anosov
diffeomorphisms are topologically mixing \cite[Corollary 18.3.5]{KH} and thus
have the specification property \cite[Theorem 18.3.9]{KH}. For the latter, 
\cite[Theorem 20.1.3]{KH} proves existence and uniqueness of the measure of
maximal entropy.  Finally, it has full support due to the estimate from below of
the measure of a dynamical ball \cite[Lemma 20.1.1]{KH}.  By uniqueness, $\mu$ 
is also $\a$ -invariant. 
Indeed, for any $b \in \Zk$ by commutativity $b_\ast \mu$ is also $a$-invariant and
has the same entropy. Hence $\mu$ lifts to an $\a'$ -invariant measure $\mu'$ on the 
suspension manifold $M'$. Clearly $\mu'$ is ergodic with respect to $\a'$ and 
has full support on $M'$.

We will now verify that the suspensions of nondegenerate TNS Cartan $\Zk$ actions
satisfy the  assumptions of Theorem \ref{technical}. First we notice the easy 

\begin{lemma}
The
suspension action of a $\Z^k$ action all of whose non-trivial elements act by Anosov
diffeomorphisms is  totally Anosov.
\end{lemma}

\begin{proof} If $a \in \Zk$ is an Anosov diffeomorphisms then any non-trivial element
on the line $\R \cdot a$  is an Anosov element for the suspension action. 
Since all non-trivial 
elements of the original action are assumed to be Anosov it is clear that $\R \cdot 
\Z ^k$ forms a dense set in $\R ^k$.   \QED
\end{proof}

Next we will verify the two  hypotheses on transitivity of Theorem \ref{technical} adapting 
an argument of \cite{KS96}. First we check the transitivity assumption that every
Lyapunov hyperplane contains a generic one-parameter subgroup with a dense orbit.

\begin{lemma}
For every Lyapunov hyperplane $\partial H$ 
almost every element in $\partial H$ is transitive.
\end{lemma}

\begin{Proof}
Fix a Lyapunov half-space $H$.
Pick a generic element $a \in \partial H$ and an Anosov element $c \in (-H)$  
so close to $a$ as to preserve the signs of all Lyapunov exponents nonzero 
on $a$. In particular we get $\ws_a=\ws_c$.  
We can take element $c$ to be a time-$t$ map of an Anosov element $b$ which
fixes the fibers of the suspension and induces an Anosov diffeomorphism
on them. Note that this Anosov diffeomorphism is transitive. Indeed, it preesrves a
finite measure of full support. Thus the nonwandering set is the whole manifold. It
is well-known that this implies topological transitivity \cite[Corollary 18.3.5]{KH}.
As above, such an element has a unique measure $\nu$ of maximal entropy
which is also $\a$ -invariant and lifts to an $\a'$ -invariant measure $\nu'$ with 
full support on the suspension manifold $M'$. We will use measure $\nu'$  for 
the rest of the proof.

Birkhoff averages with respect to $a$ of any continuous function are constant on the
leaves of $\ws _a$. Since such averages generate the algebra of $a$-invariant functions
we conclude that the partition $\xi _{a}$ into ergodic components of $a$ is coarser than
the measurable hull $\xi(\ws_a)$ of the foliation $\ws_a$ which coincides with the Pinsker 
algebra $\pi(c)$ (\cite{LY1}, Theorem B). Thus we conclude

$$\xi_a \leq \xi(\ws _a)=\xi(\ws _c)=\pi(c).$$

This shows the partition $\xi _{\partial H}$ into ergodic components
of $\partial H$ is coarser than $\xi(\ws _c) = \pi(c)$.
Since the Pinsker algebra of $b$ with respect to $\nu$ on $M$ is trivial, 
the Pinsker algebra of $c$ (with respect to $\nu'$) is coarser than the partition into 
the fibers of the suspension. Hence so is $\xi _{\partial H}$, the partition into the 
ergodic components of $\partial H$. Hence we can 
project $\xi _{\partial H}$ along the fibers of the suspension to the factor. The factor
is $\Tk$ with the standard $\rk$ action by translations. Since every element of $\Zk$
is Anosov, no element in $\Zk$ can belong to a Lyapunov hyperplane. Thus $\partial H$
contains no element in $\Zk$ and hence $\partial H$ action on $\Tk$ is uniquely ergodic. 
This implies that the projection of $\xi _{\partial H}$ to $\Tk$ is trivial, and hence so is 
$\xi _{\partial H}$ itself. This establishes the ergodicity of $\partial H$ with respect to 
$\nu'$. It is general that if an abelian group acts ergodically then so does a.e.  
element in the group. Hence we get transitivity of almost every element.
\end{Proof}
\QED

Now we verify the second transitivity assumption of Theorem \ref{technical}. 
Consider two (negative) Lyapunov half spaces $H_i$ and $H_j \not= \pm H_i$
and denote by $\w _i$ and $\w _j$ the corresponding coarse Lyapunov foliations.
By nondegeneracy, we can choose an element $a \in \partial H_i \cap \partial H_j$ 
which does not belong to any other Lyapunov hyperplane. Then we can take an 
Anosov element $c$ in  $-(H_i \cap H_j)$ so close to $a$ as to preserve the signs 
of all Lyapunov exponents nonzero on $a$. By the choice of $a$, any Lyapunov 
exponent zero on $a$ has either $\partial H_i$ or $\partial H_j$ as the kernel. 
Since the action is TNS, such an exponent must be negative on either $H_i$ or $H_j$, 
and hence is positive on $c$. Thus $\eu_c= \eu _a \oplus E_i \oplus E_j$ and 
$\ws _c =\ws_a$.  Therefore, as in the proof of the previous lemma, we get the following
inequalities 
$$\xi_a \leq \xi(\ws _a)=\xi(\ws _c)=\pi(c).$$ 
Again as in the previous lemma the Pinsker algebra of $c$ (with respect to the measure
of the maximal entropy for $c$) is coarser than the partition into the fibers of the suspension.
Thus the ergodic components of $a$ consist of the whole fibers. Since the conditional measures
on the fibers have full support, the closure of $a$-orbit of a typical point contains the whole 
fiber through that point and, in particular, $\w _{ij}$.
\QED

We conclude from Theorem \ref{Main} that $\a'$ is $\ci$ conjugate 
to an algebraic action. Now  Theorem \ref{MainZ} follows from

\begin{lemma}
If the suspension of a $\Z^k$ Anosov action is $\ci$-conjugate to an algebraic action 
then the original $\Z^k$-action is  $\ci$-conjugate to an action by automorphisms of
infranilmanifolds
\end{lemma}

\begin{proof}
The induced $\R ^k$-action  contains the original $\Z^k$-action by restricting to
a fiber. Since the algebraic $\R ^k$-action expands and contracts a smooth Riemannian
metric such that
\[||a_\ast(v)||=e^{\chi (a)} ||v|| \qquad \text{ for any } \; v \in E_\chi \]
so does the orginal $\Z^k$-action. Also notice that the non-orbit coarse Lyapunov spaces
are all tangent to the fibres of the suspension. Then  it follows by the same argument 
as at the end of the proof of Theorem 1.4 that the $\Zk$-action  is algebraic.
It is well-known that the only algebraic $\Zk$ actions are the ones by 
automorphisms of infranilmanifolds.  We refer to 
\cite[Proposition 3.13]{GS} for a proof.  
\end{proof}
\QED


\bibliographystyle{alpha}

\end{document}